\documentclass[11pt,a4paper]{amsart}
\usepackage[english]{babel}
\usepackage[utf8]{inputenc}
\usepackage[foot]{amsaddr}
\usepackage{microtype}
\usepackage{amsmath}
\usepackage{amsfonts}
\usepackage{bbm}
\usepackage{amssymb,amsthm}
\usepackage{siunitx}
\usepackage[format=plain, indention=.5cm]{caption}
\usepackage{subcaption}
\usepackage{graphicx}
\usepackage{color}
\usepackage[usenames,dvipsnames]{xcolor}
\usepackage{textcomp}
\usepackage[a4paper]{geometry} 
\usepackage[english]{babel}
\usepackage{mathtools}
\usepackage[
final,
colorlinks,
linkcolor = BrickRed,
citecolor = OliveGreen,
filecolor = black,
urlcolor = BrickRed,
bookmarks=true,
bookmarksnumbered=true,
bookmarksopen=true,
bookmarksopenlevel = 0,
hypertexnames=false,
]{hyperref}
\usepackage[
style=alphabetic,
giveninits=true,
maxbibnames=6,
minbibnames=6,
]{biblatex}

\renewbibmacro*{doi+eprint+url}{%
	\printfield{doi}%
	\newunit\newblock
	\iffieldundef{doi}
	{\usebibmacro{url+urldate}}
	{}%
	\newunit\newblock
	\usebibmacro{eprint}%
}
\usepackage[capitalise,noabbrev,nameinlink]{cleveref}
\crefname{figure}{Figure}{Figures}
\Crefname{figure}{Figure}{Figures}
\crefname{table}{Table}{Tables}
\Crefname{table}{Table}{Tables}

\usepackage{siunitx}
\newtheorem{thm}{Theorem}[section]

\usepackage{float}
\usepackage{multirow}
\usepackage{xcolor}
\usepackage{tikz}
\usetikzlibrary{shapes}
\usetikzlibrary{plotmarks}
\usetikzlibrary{calc}
\usetikzlibrary{spy}
\usetikzlibrary{fit}
\usetikzlibrary{arrows.meta}
\usetikzlibrary{snakes}
\usetikzlibrary{patterns}
\usetikzlibrary{positioning,calc}
\usetikzlibrary{decorations.markings} 
\usetikzlibrary{backgrounds}
\usetikzlibrary{shapes.misc}
\usetikzlibrary{patterns}
\usepackage{pgfplots}
\usepgfplotslibrary{colormaps}
\usepgfplotslibrary{fillbetween}

\usepgfplotslibrary{external}
\usepackage{booktabs}
\usepackage{comment}

\newcommand{\R}{\mathbb{R}}

\newcommand{\diag}{\operatorname{diag}}
\newcommand{\tr}{\operatorname{tr}}
\newcommand{\adj}{\operatorname{adj}}

\newcommand{\sign}{\operatorname{sign}}
\newcommand{\cc}{{\operatorname{c}}}
\newcommand{\pc}{{\operatorname{pc}}}

\newcommand{\Rplus}{\R_{+}}
\newcommand{\NN}{\mathcal{NN}}
\newcommand{\Loss}{\mathcal{L}}

\newcommand{\Od}{\mathcal{O}(d)}
\newcommand{\SOd}{\mathcal{SO}(d)}
\newcommand{\M}{\mathcal{M}}
\newcommand{\m}{m}

\newcommand{\Pid}{\Pi(d)}

\newcommand{\Sd}{\mathcal{S}(d)}
\newcommand{\gensym}{\mathcal{S}}

\newcommand{\sym}{{\rm sym}}
\newcommand{\ineq}{{\rm ineq}}
\newcommand{\mse}{{\rm mse}}

\definecolor{unia-purple}{HTML}{3C2673}
\definecolor{unia-pink}{HTML}{8C358B}
\definecolor{unia-yellow}{HTML}{FFDE39}
\definecolor{unia-orange}{HTML}{EF7D00}
\definecolor{unia-red}{HTML}{B93131}
\definecolor{unia-lightgreen}{HTML}{009B5D}
\definecolor{unia-green}{HTML}{13574D}
\definecolor{unia-lightblue}{HTML}{28B9D3}
\definecolor{unia-blue}{HTML}{0099DF}
\definecolor{unia-darkblue}{HTML}{3C5697}
\definecolor{unia-gray}{HTML}{b3b3b3}
\colorlet{coltrainLoss}{unia-green}
\colorlet{colvalLoss}{unia-orange}
\colorlet{colvalLossMse}{unia-lightgreen}
\colorlet{colvalLossSym}{unia-red}
\colorlet{colvalLossIneq}{unia-blue}
\colorlet{colPhi}{gray}
\colorlet{colPhiPc}{unia-lightgreen}
\colorlet{colPhiPcSVPCLP}{unia-lightgreen}
\colorlet{colPhiPcPred}{unia-pink}

\title[Neural Network Compression of Polyconvex Envelopes]{Compression of Polyconvex Envelopes of Isotropic Functions via Monotonic Input Convex Neural Networks}
\author[]{T.~Neumeier$^{*}$, J.~Salmon$^{*}$}
\address{${}^{*}$ Institute of Mathematics, University of Augsburg, Universit\"atsstr.~12a, 86159 Augsburg, Germany}
\email{\{timo.neumeier, julian.salmon\}@uni-a.de}
\date{\today}

\begin{document}
	
	\thanks{The authors gratefully acknowledge funding from the German Research Foundation (DFG) within the Priority Programme 2256 \emph{Variational Methods for Predicting Complex Phenomena in Engineering Structures and Materials} (project number 441154176, reference IDs PE1464/7-2 and PE2143/5-2). 
	Furthermore, the authors would also like to thank the Bavarian State Ministry of Science and the Arts for funding the Augsburg AI Production Network as part of the High-Tech Agenda Plus.}
	
	\begin{abstract}
		This work presents a novel neural-network compression approach for polyconvex envelopes of isotropic functions. 
		The approach relies on a classical sufficient criterion for polyconvexity and is particularly suited for the representation of determinant-constrained energy densities arising in non-linear elasticity. 
		Compared with existing compression methods based on the necessary and sufficient characterisation of polyconvex isotropic functions, the proposed framework reduces computational costs, due to the domain reduction through the restriction to the positive octant in the singed singular value space. 
		The underlying neural-network architecture employs input-convex neural networks (ICNNs) with non-negative weight constraints to enforce the required convexity and monotonicity properties. 
		The additional symmetry and inequality conditions characterising the polyconvex envelope are incorporated weakly through the loss function during training. 
		Although the employed criterion is only sufficient and thus generally yields only a lower bound on the polyconvex envelope, numerical experiments based on the classical Saint~Venant--Kirchhoff energy demonstrate that the proposed approach produces accurate approximations in practice while offering a computationally more efficient alternative to existing methods. 
	\end{abstract}
	
	\maketitle
	
	{\tiny {\bf Key words.} Polyconvexity, input convex neural network, monotonicity, relaxation}\\
	\indent
	{\tiny {\bf AMS subject classifications.} {\bf 49J45}, {\bf 49J10}, {\bf 74G65}, {\bf 74B20}, {\bf 68T07}
	}
		
\section{Introduction}
Many problems in non-linear elasticity can be formulated as variational minimisation problems of the form
\begin{equation*} 
	I(u) = \int_{\Omega} W(\nabla u)\,\mathrm{d}x,
\end{equation*}
where \(\Omega \subset \R^d\) in spatial dimension \(d \in \{2,3\}\), denotes the reference domain, \(u\colon \Omega \to \R^d\) an admissible deformation, and \(W \colon \R^{d \times d} \to \R_{\infty} \coloneqq \R \cup \{\infty\}\) the energy density function. 
In practical applications, the function value \(W(F)=\infty\) is used to model physically inadmissible deformation states, for example through orientation-preserving determinant constraints. 

Many constitutive models arising in non-linear elasticity, phase transformations, damage mechanics and fracture are inherently non-convex, see, for example, \cite{Bal76,BalJam1987,Mue99,KinJamLusEri1993,Bha2003,Ped97,BalOrt:2012:riv,Rao86,BerBoeSil07}. 
Consequently, the associated variational problems may fail to admit minimisers, while numerical approximations often exhibit mesh dependence, reduced robustness and pronounced sensitivity with respect to discretisation and material parameters.
\smallskip

A classical remedy is provided by relaxation theory. 
That is, instead of the original non-convex energy density \(W\), one considers a suitable semiconvex envelope, thereby obtaining a relaxed problem that admits minimisers and captures the effective behaviour of oscillatory minimising sequences. 
Among the various notions of semiconvexity, polyconvexity, introduced by Ball in the works \cite{Bal76,Bal77,Bal02}, constitutes a physically meaningful sufficient condition for weak lower semicontinuity of the functional and therefore for the existence of minimisers \cite{Dac08}.
Consequently, the polyconvex envelope plays a central role in the relaxation of variational problems.
	
Explicit analytical representations of semiconvex envelopes are available only in a limited number of special cases, see, for example, \cite{Koh1991,DeSDol:2002:mrn,Sil2007,KohStr86a,ConDolKre2011,ConDol2014,ConOrt2015,ConDol2017,LeDRau1994,LeRao:1995:qes,LeDRao1995,MarGhiKoeBalSanNef2025}. 
For general non-convex energy densities, however, closed-form expressions for the corresponding semiconvex envelopes are typically unavailable, necessitating computational relaxation procedures.

This has motivated the development of a broad range of computational semiconvexification methods, see, for example, \cite{Dol99,DolWal00,Bar04,BarCarHacHop04,AubFagOrt2003,ConDol18,BalKohNeuPetPet23,KNMPPB2022,KNPPB24,KNPPB2023,ObeRua:2017:pde,KumVidKoc:2020:ant}. 
In the context of polyconvexity, dedicated polyconvexification algorithms have been proposed in \cite{Bar05,EneBosGri13,BosEneGri15,FanHenKruMurWei2026}. 
For isotropic functions, dimension-reduced formulations have recently been developed in \cite{NeuPetPetWie24} based on the characterisation of polyconvexity through signed singular values established in \cite{WiePet26}.

The computation of semiconvex envelopes is computationally demanding, since the convex envelope construction is inherently non-local.
Moreover, the associated numerical procedures typically require the discretisation of the \(d \times d\)-dimensional deformation gradient space, leading to a substantial computational burden.
In the polyconvex setting, this difficulty is further aggravated by the fact that convexification is performed in a lifted space associated with the minors of the deformation gradient. 
As a result, the dimensionality of the convexification problem increases substantially, particularly in three spatial dimensions. 
Although the signed singular value characterisation considerably alleviates this issue for isotropic functions and has enabled efficient dimension-reduced polyconvexification algorithms, repeated evaluations of polyconvex envelopes remain prohibitively expensive.
The computational burden becomes particularly pronounced in applications involving concurrent relaxation within parameter-dependent boundary value problem simulations, where polyconvex envelopes must be evaluated repeatedly for varying material parameters and loading states.
\smallskip

Recent developments in machine learning have provided powerful tools for constructing structure-preserving surrogate models. 
In particular, Input Convex Neural Networks (ICNNs) \cite{AmoXuKol17} enable the incorporation of convexity directly into the network architecture and have therefore attracted significant attention in constitutive modelling and computational mechanics. 
Various neural-network-based formulations enforcing polyconvexity through sufficient criteria have been proposed, see, for example, 
\cite{KleFerMarNefWee22,KleOrtMarWee2022,LinKleKalBruWeeKae23,VijRusPauBes2025,KleRotValWee2023,ZheKocKum2024,KleHosKikKanRudGil2025}. 
Tailored architectures for isotropic hyperelasticity based on the sufficient and necessary characterisation of polyconvexity established in \cite{WiePet26} were introduced in \cite{GeuKurWieMos2025} and subsequently extended to incompressible formulations in \cite{GeuKurWieZlaCanMos2026}.

While these approaches focus on the direct learning of constitutive models, the compression of precomputed polyconvex envelopes was considered in \cite{BalNeuPetPet2025}. 
There, a property-preserving neural-network architecture was proposed for the compression of polyconvex envelopes of isotropic functions based on the sufficient and necessary characterisation from \cite{WiePet26}. 
The resulting surrogate representation acts on the signed singular value space and preserves the structural properties of the polyconvex envelope through a combination of architectural constraints and suitably designed loss terms.

The present contribution builds upon this framework and introduces an alternative compression strategy based on the classical sufficient polyconvexity criterion by Ball \cite{Bal76} for determinant-constrained isotropic functions. 
In contrast to \cite{BalNeuPetPet2025}, the proposed approach operates on the singular value space and therefore benefits from a further reduction of the characterising domain.
The corresponding neural-network architecture combines convexity-preserving structures (ICNNs) with additional monotonicity constraints inherent in Ball's criterion, similar to the monotonicity-preserving architectures considered in \cite{GeuKurWieZlaCanMos2026,KleHosKikKanRudGil2025}.
As a consequence, the resulting surrogate automatically satisfies the structural requirements of the underlying sufficient criterion for polyconvexity.
However, the additional monotonicity assumptions render the proposed representation more restrictive than the equivalent characterisation from \cite{WiePet26}. 
Consequently, the compressed representation may, in general, provide only a lower bound on the polyconvex envelope. 
Nevertheless, the reduction from the signed singular value space to the singular value space substantially decreases the computational effort.

The proposed methodology is investigated for determinant-constrained isotropic energy densities whose polyconvex envelopes are not available in analytical form. 
The numerical experiments demonstrate that the reduced representation yields substantial computational savings while maintaining the accuracy of state-of-the-art approaches based on the full signed singular value characterisation.
This makes the proposed surrogate particularly attractive for applications requiring repeated evaluations of polyconvex envelopes.
\smallskip

The remainder of this paper is organised as follows. 
\cref{sec:theory} reviews the relevant criteria for polyconvexity of isotropic functions and the associated envelope characterisations.
\cref{sec:nnarchitecture} introduces the corresponding property-preserving neural-network architectures. 
Finally, \cref{sec:numerical-experiments} presents numerical experiments and compares the proposed compression approach with existing methods.

\section{Sufficient and Necessary Criteria for Polyconvexity of Isotropic Functions} \label{sec:theory}
This section introduces the theoretical foundations of polyconvexity as well as sufficient and necessary criteria for the polyconvexity of isotopic functions. 
Let \(d \in \{2,3\}\) and let \(W \colon \R^{d \times d} \to \R_{\infty}\) denote a energy density function that maps \(d \times d\)-matrices to real scalars or infinity. 
We denote the determinant of a matrix \(F\) by \(\det(F)\) and its adjugate by \(\adj(F)\). 
Polyconvexity is characterised through the \emph{minors} of the deformation gradient, represented by the mapping \(\M \colon \R^{d \times d} \to \R^{K_d}\), where \(K_{2} = 5\) and \(k_{3} = 19\), defined by
\begin{equation*} 
	\mathcal{M}(F)=
	\begin{cases}
		(F, \, \det(F)) & \text{ if } d = 2,\\
		(F, \, \adj(F), \, \det(F)) & \text{ if } d = 3.
	\end{cases}
\end{equation*}
Specifically, a function \(V\colon\R^{d\times d}\to \R_\infty\) is said to be \emph{polyconvex} if there exists a convex function \(G \colon \R^{K_d} \to \R\) such that
\begin{equation*}
	V(F) = G(\M(F)) \qquad \text{ for all } F\in\R^{d\times d}.
\end{equation*} 

Let \(\Od\) denote the orthogonal group, \(\SOd\) the special orthogonal group and \(\Sd \subseteq \{0,1\}^{d\times d}\) the group of permutation matrices.
Throughout this work, we focus on \emph{isotropic} energy densities \(W\colon\R^{d\times d}\to \R_\infty\), i.e.~for all \(F\in\R^{d\times d}\) and for all \(R_1, R_2\in \SOd\) it holds 
\begin{equation*}
	W(F) = W(R_1\, F\, R_2).
\end{equation*}
Note that this definition of isotropy incorporates objectivity and full material symmetry, following the convention outlined in \cite{Bal76}, and is also referred to as \(\SOd\times\SOd\)-invariance. 

A classical result, see e.g.~\cite[Proposition~5.31]{Dac08}, states that isotropic functions admit a dimension-reduced representation in terms of signed singular values.
Based on the signed singular value decomposition, i.e.~for all \(F\in\R^{d\times d}\) there exist \(R_1, R_2\in \SOd\) and \(\hat{\nu}\in \R^d\) such that
\begin{equation*}
	F = R_1 \, \mathrm{diag}(\hat{\nu}) \, R_2, 
\end{equation*}
it is possible to derive the dimension reduced representation.
In general, the signed singular values \(\hat{\nu}\in \R^d\) are only unique up to transformations included in the symmetry group consisting of orientation preserving signed permutations, formally described by
\begin{equation*}
	\Pid = \left\{\diag(\varepsilon )\, S \in \Od \;\middle|\; S \in \Sd,\, \varepsilon \in \{ \pm 1\}^d,\, \varepsilon_{1} \cdots \varepsilon_{d} = 1\right\}.
\end{equation*}
A canonical representative might be obtained through the signed singular value mapping.
We introduce the singular value mapping \(\sigma\colon\R^{d\times d}\to\R^d_+\), where \(R_{+} = \{x \geq 0\}\), defined by
\begin{equation*}
	\sigma(F) = \left[\sigma_1(F), \, \ldots, \, \sigma_d(F)\right]^\top,
\end{equation*}
with \(0\leq \sigma_1(F) \leq \ldots \leq \sigma_d(F)\). 
The signed singular value mapping \(\nu\colon\R^{d\times d}\to\R^d\) is then defined by
\begin{equation*}
	\nu(F) = \left[\sign(\det(F)) \, \sigma_1(F), \, \sigma_2(F), \, \ldots, \, \sigma_d(F)\right]^\top.
\end{equation*}
Every isotropic function \(W\colon\R^{d\times d}\to\R_\infty\) is in one-to-one correspondence with a \(\Pi_d\)-invariant function \(\Phi\colon\R^d\to\R_\infty\), i.e.~\(\Phi(\hat{\nu}) = \Phi(S \, \hat{\nu})\) for all \(S\in\Pid\), through the identifications 
\begin{equation*}
	W(F) = \Phi(\nu(F)) \quad \text{ for all } F\in\R^{d\times d} \quad \text{ and vice versa } \quad \Phi(\hat{\nu}) = W(\diag(\hat{\nu})) \quad \text{ for all } \hat{\nu} \in \R^{d}.
\end{equation*}

Following \cite{WiePet26}, polyconvexity can equivalently be transferred to functions \(\Phi\colon\R^d\to\R_\infty\) acting on the signed singular values.
To this end, we introduce the vector-valued minor mapping \(\m\colon\R^d\to\R^{k_d}\), where \(k_2=3\) and \(k_3=7\), defined by
\begin{equation*}
	\m(\hat{\gamma}) = 
	\begin{cases}
		(\hat{\gamma}_1, \, \hat{\gamma}_2, \, \hat{\gamma}_1 \, \hat{\gamma}_2) & \text{ if } d = 2, \\
		(\hat{\gamma}_1, \, \hat{\gamma}_2, \, \hat{\gamma}_3,\, \hat{\gamma}_2 \, \hat{\gamma}_3, \, \hat{\gamma}_3 \, \hat{\gamma}_1, \, \hat{\gamma}_1 \, \hat{\gamma}_2, \, \hat{\gamma}_1 \, \hat{\gamma}_2 \, \hat{\gamma}_3) & \text{ if } d = 3.
	\end{cases}
\end{equation*}
A function \(\Phi\colon\R^d\to\R_\infty\) is said to be \emph{polyconvex} if there exists a convex function \(g \colon \R^{k_d} \to \R_{\infty}\) such that
\begin{equation*}
	\Phi(\hat{\nu}) = g(\m(\hat{\nu})) \qquad \text{ for all } \hat{\nu}\in \R^d.
\end{equation*}

The following theorem provides a necessary and sufficient criterion for the polyconvexity of isotropic functions and establishes the equivalence of polyconvexity in the deformation gradient and signed singular value representations.
The theorem below was stated by Wiedemann and Peter in \cite{WiePet26}, where lower semicontinuity is incorporated explicitly into the characterisation.
Since lower semicontinuity is automatically satisfied within our neural-network framework considered here, we omit this aspect and refer to \cite{WiePet26} for the details of the statement.
\begin{thm}[{\cite[Theorem~1.4]{WiePet26}}] \label{thm:isocrit}
	Let \(W\colon\R^{d\times d}\to \R_\infty\) be isotropic and \(\Phi\colon\R^{d}\to\R_\infty\) be the associated \(\Pid\)-invariant function.
	Then the following are equivalent:
	\begin{itemize}
		\item \(W\) is polyconvex,
		\item \(\Phi\) is polyconvex,
		\item there exists a convex function \(g \colon \R^{k_d} \to \R\) such that
		\begin{equation}
			W(F) = \Phi(\nu(F)) = g(\m({\nu(F)}))
		\end{equation}
		for all \(F\in\R^{d\times d}\) and the composition \(g \circ m\) is \(\Pid\)-invariant.
	\end{itemize}
\end{thm}

While \cref{thm:isocrit} provides an equivalent characterisation of polyconvex isotropic functions, many explicit constructions in non-linear elasticity are based on the classical sufficient criterion due to Ball \cite{Bal76}, which is formulated directly in terms of the singular values.
\begin{thm}[{\cite[Theorem~5.2]{Bal76}}] \label{thm:ballcrit}
	Let \(U = \{F \in \R^{d \times d} \mid \det (F) \in \Rplus\}\) and let \(\bar{g} \colon \R^{k_d}_+ \to \R\) be a convex function which is non-decreasing in the first \(k_d-1\) arguments.
	Assume that the composition \(\Upsilon = \bar{g} \circ m\) is \(\Sd\)-invariant and satisfies
	\begin{equation*}
		W(F) = \Upsilon(\sigma(F)) = \bar{g}(\m(\sigma(F)))
	\end{equation*}
	for all \(F\in U\). 
	Then the function \(W\) is polyconvex on \(U\). 
\end{thm}
Moreover, under the assumptions of \cref{thm:ballcrit}, the result extends naturally to extended-real-valued functions incorporating orientation-preserving determinant constraints. 
In particular, let \(\Upsilon \colon \R_{+}^{d} \to \R_{\infty}\) be of the form \(\Upsilon(\hat{\sigma}) = \bar{g}(\m(\hat{\sigma}))\) for all \(\hat{\sigma} \in \R^d_{+}\) where \(\bar{g}\) satisfies the assumptions from \cref{thm:ballcrit}. 
Then the function
\begin{equation} \label{eq:ballcritextended}
	W(F) = 
	\begin{cases}
		\Upsilon(\sigma(F)) & \text{ if } \det(F) > 0, \\
		\infty & \text{ otherwise}
	\end{cases}
\end{equation}
is polyconvex on \(\R^{d\times d}\). 	
\smallskip

While \cref{thm:ballcrit} provides a convenient sufficient condition of polyconvex functions, a dimension reduced computation of polyconvex envelopes relies on the characterisation via \cref{thm:isocrit}. 
In particular, the envelope construction can be reduced to the convexification of a lifted function in the signed singular value space.
To this end, we briefly recall the definition of the polyconvex envelope.

The polyconvex envelope of the energy density \(W\) is the largest polyconvex function below \(W\), i.e.~the mapping \(W^{\pc} \colon \R^{d\times d} \to \R_\infty\) defined by
\begin{equation*}
	W^{\pc}(F) = \sup \left\{ V(F) \;\middle|\;  V \colon \R^{d \times d} \to \R_\infty,\ V \text{ polyconvex, } V \leq W \right\}.
\end{equation*}
Correspondingly, for a function \(\Phi \colon \R^d \to \R_\infty\), the polyconvex envelope is defined as the largest polyconvex function below \(\Phi\), namely
\begin{equation} \label{eq:defPhipcenvelope}
	\Phi^{\pc}(\hat{\nu}) = \sup \left\{ \Psi(\hat{\nu}) \;\middle|\;  \Psi \colon \R^{d} \to \R_\infty,\ \Psi \text{ polyconvex, } \Psi \leq \Phi \right\}.
\end{equation}

For isotropic functions \(W\) and the associated \(\Pid\)-invariant functions \(\Phi\), \cite[Proposition~2.4]{NeuPetPetWie24} shows that the polyconvex envelopes can be determined through the convexification of the function \(h \colon \R^{k_d} \to \R_{\infty}\) defined on the lifted signed singular values by
\begin{equation} \label{eq:h}
	h(x) = 
	\begin{cases}
		W(F) &  \text{ if } x = \m(\nu(F)), \\
		\infty & \text{ otherwise}.
	\end{cases}
\end{equation}
In particular, for all \(F\in\R^{d\times d}\) the natural identification
\begin{equation} \label{eq:WpcPhipchc}
	W^\pc(F)=\Phi^\pc(\nu(F))=h^\mathrm{c}(\m(\nu(F)))
\end{equation}
holds.
\medskip

The characterisation \eqref{eq:WpcPhipchc} provides a viable strategy for determining the polyconvex envelope via a standard convexification in the lifted signed singular value space.
In the following, we consider an alternative representation motivated by the sufficient criterion of \cref{thm:ballcrit}. 
Restricting attention to determinant-constrained isotropic functions \(W\colon\R^{d\times d}\to \R_\infty\), we assume that the polyconvex envelope of \(W\) admits a representation of the form
\begin{equation} \label{eq:Wpc=Upsilonpc}
	W^{\pc}(F) 
	= 
	\begin{cases}
		\Upsilon^{\pc}(\sigma(F)) & \text{ if } \det(F) > 0, \\
		\infty & \text{ otherwise}, 
	\end{cases}
\end{equation}
where \(\Upsilon^{\pc}\) is \(\Sd\)-invariant and satisfies
\begin{equation} \label{eq:Upsilonpc=barhc}
	\Upsilon^{\pc}(\hat{\sigma}) = \bar{h}^{\cc}(\m(\hat{\sigma})), 
\end{equation}
with \(\bar{h}^{\cc}\colon \R^d_{+} \to \R_{\infty}\) convex and non-decreasing in the first \(k_d-1\) arguments.
It should be stressed, that technically \(\Upsilon^{\pc}\) is \emph{not a polyconvexification}. 
Rather, it serves as the singular value representation of the polyconvex envelope and is a notation used to distinguish it from \(\Phi^{\pc}\).
Indeed, under the assumption \eqref{eq:Wpc=Upsilonpc}, the singed singular value polyconvex envelope is related by
\begin{equation*}
	\Phi^{\pc} (\hat{\nu}) = 
	\begin{cases}
		\Upsilon^{\pc}(\sigma(\diag(\hat{\nu}))) & \text{ if } \nu_{1} \cdot \ldots \cdot \nu_{d} > 0, \\
		\infty & \text{ otherwise}.	
	\end{cases}	
\end{equation*}
The composition with \(\sigma \circ \diag\) ensures that \(\Upsilon^{\pc}\) gets the singular values as argument, i.e.~the input vector \(\hat{\nu}\) is mapped to the positive octant as it was the convention of the singular value mapping \(\sigma\). 
Consequently, \(\Upsilon^{\pc}\) may be viewed as the restriction of \(\Phi^{\pc}\) to the singular value cone.

It is important to note that the representation of the form \eqref{eq:Wpc=Upsilonpc} is an assumption. 
Indeed, the monotonicity requirements imposed on \(\bar{h}^{\cc}\) are not implied by the characterisation \eqref{eq:WpcPhipchc}. 
If these conditions are not satisfied, the resulting function obtained from the monotone convexification, i.e.~a monotonous convex envelope of the function \(h\) from \eqref{eq:h} over \(\R^{k_d}_{+}\), the formulation \(\Upsilon^{\pc}\) in \eqref{eq:Wpc=Upsilonpc} still characterises a polyconvex function below \(W\), but in general constitutes only a lower bound of the polyconvex envelope \(W^{\pc}\).

The representation of the polyconvex envelope of \(W\) via \(\Upsilon^{\pc}\) as in \eqref{eq:Wpc=Upsilonpc} provides the foundation for the alternative compression approach considered in this work. 
While the additional monotonicity requirements may, in general, render the resulting approximation a lower bound on the polyconvex envelope, the associated reduction from the signed singular value space to the singular value space yields a considerable computational advantage and motivates the neural-network architecture introduced in the following section.

\section{Property-Preserving Neural Network Architectures} \label{sec:nnarchitecture}
This section introduces the neural-network-based representations of polyconvex envelopes of isotropic functions based on the characterisations presented in the previous section. 
In addition to the compression of the signed singular value polyconvex envelope \(\Phi^{\pc}\) from \eqref{eq:WpcPhipchc}, as proposed in \cite{BalNeuPetPet2025}, we consider the singular value based representation \(\Upsilon^{\pc}\) from \eqref{eq:Wpc=Upsilonpc}. 
The resulting neural network approximations are denoted by \(\Phi^{\pc}_{\NN}\) and \(\Upsilon^{\pc}_{\NN}\), respectively.

The approximation \(\Phi^{\pc}_{\NN}\) is based on the necessary and sufficient characterisation of \cref{thm:isocrit}, whereas \(\Upsilon^{\pc}_{\NN}\) follows the representation associated to \cref{thm:ballcrit}, i.e.~the Ball-criterion-based characterisation. 
Although the latter may, in general, provide only a lower bound on the polyconvex envelope, it benefits from a reduced characterising domain, acting solely on the
positive octant in the singed singular value space, the singular value cone.
As a consequence, substantially lower computational costs can be achieved while preserving the structural properties required by the corresponding polyconvexity criterion.

Both approximations are realised by the same underlying architecture, so called Input Convex Neural Networks, as proposed in \cite{AmoXuKol17}. 
Following this contribution, convexity is enforced through non-negativity constraints on selected weights together with convex and monotone activation functions. 
In the present setting for polyconvexity, convexity is required with respect to the minors, i.e.~the lifted (signed) singular values.

To describe both variants simultaneously, we denote by \(\hat{\Gamma}^{\pc}_{\NN} \in \{\Upsilon^{\pc}_{\NN}, \Phi^{\pc}_{\NN}\}\) the considered approximations in an abstract sense.
Furthermore, let \(\hat{\gamma}\in\R^{d}\) denote the corresponding input variable, i.e.~\(\hat{\nu} \in \R^{d}\) in the case of \(\Phi^{\pc}_{\NN}\) and \(\hat{\sigma}\in\R^d_+\) for \(\Upsilon^{\pc}_{\NN}\), and define the lifted input of the network by \(\hat{m}= \m(\hat{\gamma})\).
The application of the minors mapping \(\m\) is interpreted as a hard coded feature-extraction layer. 
For an \(L\)-layer network, the hidden activations \(z_{\ell}\) are defined recursively by
\begin{equation*}
	z_{\ell+1} = \rho_\ell \left(W_\ell^{(z)} z_{\ell} + W_\ell^{(\hat{m})} \hat{m} + b_\ell \right) 
\end{equation*}
for \(\ell = 1, \ldots, L - 1\) and the output \(z_{\ell} \in \R\).
The weights \(W_\ell^{(z)}\), the passthrough layer weights \(W_\ell^{(\hat{m})}\) and biases \(b_\ell\) are collected in the trainable parameter vector \(\theta\), with the typical convention \(z_0 = 0\) and \(W_0^{(z)} \equiv 0\).
The overall network output is denoted by
\begin{equation*} 
	\Upsilon^{\pc}_{\NN}(\hat{\sigma}; \theta) = \bar{h}^{\cc}_{\NN}(\m(\hat{\sigma});\theta) = z_{L}
\end{equation*}
and 
\begin{equation*} 
	\Phi^{\pc}_{\NN}(\hat{\nu}; \theta) = h^{\cc}_{\NN}(\m(\hat{\nu});\theta) = z_{L},
\end{equation*}
respectively. 

Note that neither \(\Phi_{\NN}^{\pc}\) nor \(\Upsilon_{\NN}^{\pc}\) explicitly encode the determinant constraint included in the function \(W^{\pc}\) as in \eqref{eq:Wpc=Upsilonpc}. 
Both networks are finite-valued representations of the corresponding characterising functions and the extension by \(\infty\) for \(\det(F)\leq 0\) should be enforced separately during the evaluation of \(W_{\NN}^{\pc}\), in agreement with \eqref{eq:Wpc=Upsilonpc}.

The convexity and monotonicity properties of \(h^{\cc}_{\NN}\) and \(\bar{h}^{\cc}_{\NN}\) required by \cref{thm:isocrit} and \cref{thm:ballcrit} are enforced through suitable restrictions on the network parameters and activations. 
According to \cite{BalNeuPetPet2025} and \cite{AmoXuKol17}, convexity of the network output with respect to the lifted (signed) singular values input arguments is guaranteed by employing convex and non-decreasing activation functions \(\rho_\ell\) and by constraining the weights in the \(z\)-path to be non-negative, i.e.~\(W_\ell^{(z)} \geq 0\) for \(\ell = 1,\ldots,L-1\).
Under these assumptions, the functions \(h^{\cc}_{\NN}\) and \(\bar{h}^{\cc}_{\NN}\) are convex with respect to their inputs.

For the representation \(\Upsilon^{\pc}_{\NN}\), the monotonicity assumptions of \cref{thm:ballcrit} must additionally be satisfied. 
Since the activation functions are already chosen to be non-decreasing, the additional monotonicity of \(\bar{h}^{\cc}_{\NN}\) in the first \(k_d-1\) arguments is obtained by enforcing the corresponding columns of the weights in the \(\hat{m}\)-path to be non-negative, i.e.~\([W_\ell^{(\hat m)}]_{:,1:k_d-1} \geq 0,\) for \(\ell = 0,\ldots,L-1\).

Consequently, \(\Phi^{\pc}_{\NN}\) satisfies the convexity requirements of \cref{thm:isocrit}, whereas \(\Upsilon^{\pc}_{\NN}\) additionally satisfies the monotonicity assumptions required by \cref{thm:ballcrit}. 
Therefore, both architectures preserve the structural properties underlying the corresponding polyconvexity criteria by construction.
\smallskip

Following \cite{BalNeuPetPet2025}, the remaining characteristic properties, namely the symmetry conditions and the envelope inequality, are incorporated into the neural-network models in a weak sense through the loss function. 
Let \(\gensym\) denote the associated symmetry group to the considered neural-network representation \(\hat{\Gamma}^{\pc}_{\NN}\), i.e.~\(\Pid\) in the case of \(\Phi^{\pc}_{\NN}\) and \(\Sd\) for \(\Upsilon^{\pc}_{\NN}\).
Furthermore, let \((\hat{\gamma}_{i}, \Phi^{\pc}(\hat{\gamma}_{i}), \Phi(\hat{\gamma}_{i}))\), for \(i = 1, \ldots, N\) denote the training data points. 
The symmetry property is enforced through the penalty term
\begin{equation*} 
	\Loss_{\sym}(\theta; \gensym) = \frac{1}{N} \sum_{i = 1}^{N} \frac{1}{\lvert \gensym \rvert} \sum_{S \in \gensym} \left(
	\hat{\Gamma}^{\pc}_{\NN}({\hat{\gamma}_{i}};\theta) - \hat{\Gamma}^{\pc}_{\NN}(S \, \hat{\gamma}_i;\theta)\right)^2.
\end{equation*}	
Moreover, the inequality property in the definition of the polyconvex envelope \(\Phi^{\pc} \leq \Phi\), is incorporated through the penalty
\begin{equation*}
	\Loss_{\ineq}(\theta) = \frac{1}{N} \sum_{i = 1}^{N}\max \bigl\{\hat{\Gamma}^{\pc}_{\NN}(\hat{\gamma}_{i};\theta)- \Phi(\hat{\gamma}_{i}),0\bigr\}^{2}.
\end{equation*} 
Combining these contributions with the mean squared approximation error for the target envelope function values yields the overall loss function
\begin{equation*}
	\Loss(\theta) = \Loss_{\mse}(\theta) + \lambda_{\sym} \, \Loss_{\sym}(\theta) + \lambda_{\ineq} \, \Loss_{\ineq}(\theta),
\end{equation*}	
where \(\lambda_{\sym}\) and \(\lambda_{\ineq}\) denote scalar penalty parameters associated with the symmetry and inequality contributions, respectively.

\section{Numerical Experiments} \label{sec:numerical-experiments}
In the following section, we focus on a well-established benchmark problem to assess the performance of the proposed compression approach and compare it with the state-of-the-art ansatz introduced in \cite{BalNeuPetPet2025}.
We consider the Saint~Venant--Kirchhoff energy density in three spatial dimensions, incorporating determinant constraints that penalise self-intersection and self-interpenetration. 
This energy density \(W\colon \R^{3 \times 3} \to \R_{\infty}\) is defined by
\begin{equation*} 
	W(F) = 
	\begin{cases}
		\frac{\mu}{4} \lvert F^\top F - \mathbb{I}_{3} \rvert^2 + \frac{\lambda}{8} \left(\lvert F\rvert^2 - 3 \right)^2 & \text{ if } \det(F) > 0, \\
		\infty & \text{ otherwise},
	\end{cases}
\end{equation*}
where \(\lvert A \rvert = \sqrt{\tr(A^\top A)}\) denotes the Frobenius norm and \(\mathbb{I}_{3} \in \R^{3 \times 3}\) the identity matrix. 
Exploiting the \(\SOd\times\SOd\)-invariance of \(W\), the function can be rephrased in the dimension-reduced representation \(\Phi \colon \R^{3} \to \R_{\infty}\), acting on the signed singular values \(\hat{\nu} = [\nu_1, \ldots, \nu_d]^\top \in \R^3\), given by 
\begin{equation*} 
	\Phi(\hat{\nu}) = 
	\begin{cases}
		\frac{\mu}{4} \, \sum_{i = 1}^{3} \left(\nu_i^2 - 1\right)^2 + \frac{\lambda}{8} \left(\lvert\hat{\nu}\rvert^2 - 3\right)^2 & \text{ if }  \nu_{1} \, \nu_{2} \, \nu_{3} > 0, \\
		\infty & \text{ otherwise}.
	\end{cases}
\end{equation*}
For the unconstrained model, the polyconvex envelope is known in closed form, see \cite{LeRao:1995:qes}. 
However, no analytical representation is known for the determinant-constrained setting considered here, making a numerical approximation of the envelope indispensable.

Owing to the isotropic structure of the energy density, the polyconvex envelope can be approximated efficiently by the signed singular value polyconvexification algorithm based on linear programming (SVPC~LP) proposed in \cite{NeuPetPetWie24}. 
The method performs a convexification of a discrete representation of the function \(h\) from \eqref{eq:h} on the minors manifold using a lifted signed singular value discretisation, i.e.~
\begin{equation*}
	W^{\pc}_{\delta}(F) = \Phi^{\pc}_{\delta}(\nu(F)) = 
	\begin{cases}
		h^\cc_{\delta}(\m(\nu(F))) & \text{ if } \det(F) > 0, \\
		\infty & \text{ otherwise}.
	\end{cases}
\end{equation*}
The resulting numerical approximation of the singed singular value polyconvex envelope is denoted by \(\Phi^{\pc}_{\delta}\).

\begin{figure}
	\centering
	\begin{tikzpicture}
	\pgfmathsetmacro{\plotskipx}{1.4}
	\pgfmathsetmacro{\plotskipy}{-7.0}
	\pgfmathsetmacro{\widthscale}{0.5}
	
	\begin{axis}[
		name=plot11,
		width=\widthscale\textwidth,
		grid=both, 
		view={35}{30},
		xmajorgrids, ymajorgrids,
		minor grid style={black!10}, major grid style={black!40},
		axis background/.style={fill=white},
		xlabel={\(\nu_1\)},	ylabel={\(\nu_2\)},
		xtick={{0.5, 1, 1.5, 2}}, 
		ytick={{0.5, 1, 1.5, 2}}, 
		zmin=-0.02, zmax=0.41,
		minor tick num=1,
		colormap/viridis, 
		mesh/ordering=x varies,
		unbounded coords=jump,
		]
		\addplot3[
		surf, 
		] 
		file {figures/data/PhisliceSTVK_3d.dat};
	\end{axis}
	\node[
	draw, 
	fill=white!30, 
	draw opacity=0.8, 
	fill opacity=0.7, 
	text opacity=1, 
	inner sep=3pt, 
	anchor=south,
	xshift=0.8cm, 
	yshift=-1.4cm, 
	] at (plot11.north) {\(\Phi\)};
	
	\begin{axis}[
		name=plot12,
		at={(plot11.south east)},
		xshift=\plotskipx cm,
		width=\widthscale\textwidth,
		grid=both, 
		view={35}{30},
		xmajorgrids, ymajorgrids,
		minor grid style={black!10}, major grid style={black!40},
		axis background/.style={fill=white},
		xlabel={\(\nu_1\)},	ylabel={\(\nu_2\)},
		xtick={{0.5, 1, 1.5, 2}}, 
		ytick={{0.5, 1, 1.5, 2}}, 
		zmin=-0.02, zmax=0.41,
		minor tick num=1,
		colormap/viridis, 
		mesh/ordering=x varies,
		unbounded coords=jump,
		]
		\addplot3[
		surf, 
		] 
		file {figures/data/PhipcsliceSTVK_3d.dat};
		
	\end{axis}
	\node[
	draw, 
	fill=white!30, 
	draw opacity=0.8, 
	fill opacity=0.7, 
	text opacity=1, 
	inner sep=3pt, 
	anchor=south,
	xshift=0.8cm, 
	yshift=-1.4cm,
	] at (plot12.north) {\(\Phi^{\pc}_{\delta}\)};
	
	\begin{axis}[
		name=plot21,
		at={(plot11.south west)},
		yshift=\plotskipy cm,
		width=\widthscale\textwidth,
		grid=both, 
		view={35}{30},
		xmajorgrids, ymajorgrids,
		minor grid style={black!10}, major grid style={black!40},
		axis background/.style={fill=white},
		xlabel={\(\nu_1\)},	ylabel={\(\nu_2\)},
		xtick={{0.5, 1, 1.5, 2}}, 
		ytick={{0.5, 1, 1.5, 2}}, 
		zmin=-0.02, zmax=0.41,
		minor tick num=1,
		colormap/viridis, 
		mesh/ordering=x varies,
		unbounded coords=jump,
		]
		\addplot3[
		surf, 
		] 
		file {figures/data/PhipcNNsliceSTVK_3d_Ball.dat};		
	\end{axis}
	\node[
	draw, 
	fill=white!30, 
	draw opacity=0.8, 
	fill opacity=0.7, 
	text opacity=1, 
	inner sep=3pt, 
	anchor=south,
	xshift=0.8cm, 
	yshift=-1.4cm, 
	] at (plot21.north) {\(\Upsilon^{\pc}_{\NN}\)};
	
	\begin{axis}[
		name=plot22,
		at={(plot21.south east)},
		xshift=\plotskipx cm,
		width=\widthscale\textwidth,
		grid=both, 
		view={35}{30},
		xmajorgrids, ymajorgrids,
		minor grid style={black!10}, major grid style={black!40},
		axis background/.style={fill=white},
		xlabel={\(\nu_1\)},	ylabel={\(\nu_2\)},
		xtick={{0.5, 1, 1.5, 2}}, 
		ytick={{0.5, 1, 1.5, 2}}, 
		zmin=-0.02, zmax=0.41,
		minor tick num=1,
		colormap/viridis, 
		mesh/ordering=x varies,
		unbounded coords=jump,
		]
		\addplot3[
		surf, 
		] 
		file {figures/data/PhipcNNsliceSTVK_3d_Iso.dat};
	\end{axis}
	\node[
	draw, 
	fill=white!30, 
	draw opacity=0.8, 
	fill opacity=0.7, 
	text opacity=1, 
	inner sep=3pt, 
	anchor=south,
	xshift=0.8cm, 
	yshift=-1.4cm, 
	] at (plot22.north) {\(\Phi^{\pc}_{\NN}\)};
		
\end{tikzpicture}
	\caption{
		Non-convex function \(\Phi\), computationally approximated polyconvex envelope \(\Phi^{\pc}_{\delta}\) and neural network compressions via \(\Upsilon^{\pc}_{\NN}\) and \(\Phi^{\pc}_{\NN}\) for the three dimensional Saint~Venant--Kirchhoff example. 
		The illustration on the cross section \((\nu_1, \nu_2, 1)\) is restricted to the positive octant of the signed singular value space.
	}
	\label{fig:STVK3Dfunctions}
\end{figure}
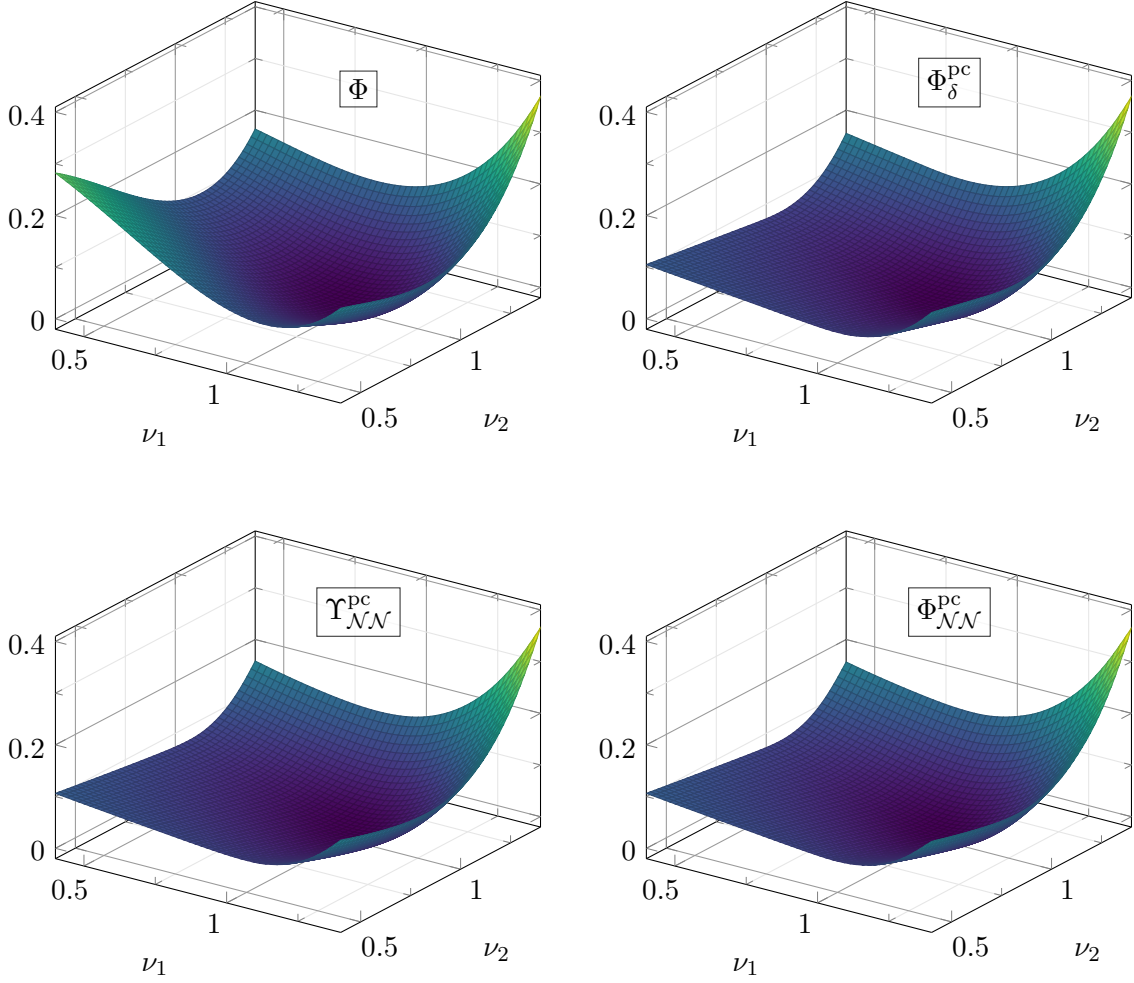

The original energy density \(\Phi\) and its numerical polyconvex envelope approximation \(\Phi^{\pc}_{\delta}\) are depicted in \cref{fig:STVK3Dfunctions}.
Since no analytical envelope is available, this benchmark provides a particularly suitable test case for the proposed compression framework: the numerically computed envelope point values obtained by SVPC~LP serve as target values that are subsequently compressed using the structure-preserving neural network approaches introduced before.

Studying the numerical envelope approximation \(\Phi^{\pc}_{\delta}\) suggests that the polyconvex envelope admits a representation of the form \eqref{eq:ballcritextended}. 
This motivates the application of the Ball-criterion-based compression approach, which amounts to learning the associated convex function \(\bar{h}^{\cc}\) from \eqref{eq:Upsilonpc=barhc} on \(\R_{+}^{3}\).

Since Ball's criterion, cf.~\cref{thm:ballcrit}, provides only a sufficient condition for polyconvexity, the additional monotonicity constraints imposed on \(\bar{h}^{\cc}\), and in the neural network may, in principle, be overly restrictive. 
Consequently, the corresponding neural network ansatz may represent only a lower bound of the true polyconvex envelope. 
For the present benchmark, however, the numerical reference approximation \(\Phi^{\pc}_{\delta}\) already satisfies the required monotonicity properties and the Ball-criterion-based neural network can be employed without introducing any observable additional approximation error in the experiments presented below.
\smallskip

The numerical setup is as follows. 
All experiments are conducted in the three-dimensional setting with Lamé parameters \(\mu = 0.4\) and \(\lambda = 0.4\).
We consider the proposed compression model \(\Upsilon^{\pc}_{\NN}\) from \eqref{eq:Upsilonpc=barhc} and the state--of--the--art model \(\Phi^{\pc}_{\NN}\) from \eqref{eq:WpcPhipchc} for comparison reasons. 
Both models are implemented in \textit{Python} using the \textit{PyTorch} framework.
The architectures for both approaches share the same principal structure. 
After the feature extraction via the hard coded minors layer, the neural network consists of a input layer of dimension \(k_{3} = 7\), corresponding to the lifted signed singular value dimension, followed by three hidden layers with \(10\), \(10\) and \(20\) neurons, respectively, and a scalar (finite) valued output layer. 

As activation function, we employ the Softplus function with parameter \(\beta=20\), which is convex and non-decreasing and therefore compatible with the theoretical requirements to ensure convexity. 
Furthermore, the weights satisfy the structure-preserving architecture constraints introduced in the previous section by projection of the relevant weights onto the positive half space via the Softplus function.
With this configuration, each network comprises \num{648} trainable parameters and apart from the additional monotonicity constraints required for the Ball-criterion-based model \(\Upsilon^{\pc}_{\NN}\), both architectures employ essentially the same implementation, ensuring a fair comparison of the respective compression strategies.

Training is performed using the Adam optimiser with learning rate \(\eta=\num{5e-3}\) and batch size \(\num{256}\). 
Early stopping with a patience of \(10\) epochs is employed, while the maximal number of training epochs is set to \num{250}. 
The penalty parameters in the loss function \(\Loss\) are chosen as \(\lambda_{\sym}=2.0\) and \(\lambda_{\ineq}=5.0\).
\smallskip

We aim to compress the polyconvex envelope on a bounded subset of the signed singular value space rather than along prescribed deformation paths. 
Consequently, the learning domain contains the signed singular values corresponding to arbitrary deformation states and is not restricted to specific deformation patterns.
The discretisation is based on the box \([\nu_{\min}, \nu_{\max}]^{3} \subset \R_{+}^{3}\) with \(\nu_{\min}=0.4\) and \(\nu_{\max}=1.4\). 
Each coordinate direction is discretised by \(75\) equidistantly distributed points.

For the Ball-criterion-based architecture \(\Upsilon^{\pc}_{\NN}\), the learning domain is given by \([\nu_{\min}, \nu_{\max}]^{3} \subset \R_{+}^{3}\), resulting in a total of \num{421875} lattice points.
In contrast, the architecture \(\Phi^{\pc}_{\NN}\) operates on the entire signed singular value space, to be more precise, the part of signed singular value space which corresponds to \(\det F > 0\), i.e.~\(\{\hat{\nu} \in \R^{d} \mid \nu_1 \, \nu_2 \, \nu_3 > 0\}\).
Therefore, all sign-preserving permutations must additionally be considered, i.e.~\(\Pid \, [\nu_{\min}, \nu_{\max}]^{3}\), yielding a total of \num{1687500} data points in \([\nu_{\min}, \nu_{\max}]^{3} \, \subset \R^{3}\).
This extension is necessary in order to satisfy the appropriate notion of convexity on the full signed singular value space.

It should be noted that throughout this work, we restrict ourselves to the finite-valued regime of the energy density and do not explicitly consider the extended-real-valued setting. 
In a practical simulation framework, the determinant constraint could instead be enforced by a dedicated preprocessing layer that detects inadmissible deformation states and returns the penalty value \(\infty\).

The restriction of \(\Upsilon^{\pc}_{\NN}\) to the positive quadrant reduces the size of the learning domain by a factor of \(4\), clearly illustrating the domain reduction achieved by the Ball-criterion-based formulation. 
This reduction substantially decreases the number of data points involved in the training process. 
The trade-off is the additional monotonicity requirement, which may in general be overly restrictive and consequently yield only a lower-bound approximation of the polyconvex envelope.

The learning data is given by tuples of the form \((\hat{\nu}, \Phi^{\pc}_{\delta}(\hat{\nu}), \Phi(\hat{\nu}))\), where the envelope values are obtained by the SVPC~LP algorithm. 
The original function values \(\Phi\) are included in order to enforce the inequality constraint appearing in the loss function. 
The resulting data set is split randomly into a training set (70\%) and a validation set (30\%). 
To reduce the influence of the random data partitions and random weight initialisations, all experiments are repeated for ten independent network realisations.

\bigskip
\cref{fig:STVK3Dfunctions} depicts the resulting averaged neural network compressions \(\Upsilon^{\pc}_{\NN}\) and \(\Phi^{\pc}_{\NN}\) on a representative cross-sectional slice in the positive octant. 
It should be noted that this domain constitutes the natural domain for \(\Upsilon^{\pc}_{\NN}\), whereas for \(\Phi^{\pc}_{\NN}\) it merely represents a restriction of the full signed singular value space used during training.
For illustration and comparison reasons, this restriction is used throughout the visualisations below.

Visually, both neural network approximations are nearly indistinguishable and closely match the reference envelope \(\Phi^{\pc}_{\delta}\).
This observation suggests that, for the determinant-constrained Saint~Venant--Kirchhoff model considered here, the Ball-criterion-based ansatz is indeed able to recover the polyconvex envelope. 
In particular, although \(\Upsilon^{\pc}_{\NN}\) is, in general, expected to provide only a lower-bound, no discrepancy between the corresponding compression and the reference envelope is observed.

\begin{table}[h]
	\centering
	\setlength{\arrayrulewidth}{0.8pt}
	\begin{tabular}{@{}p{4cm}rr}
		& \(\Upsilon^{\pc}_{\NN}\) & \(\Phi^{\pc}_{\NN}\)\\
		\midrule
		\(\Loss\) training data & \((1.1 \pm 0.6) \times 10^{-5}\) & \((1.3 \pm 0.1) \times 10^{-5}\) \\
		\(\Loss\) validation data & \((2.0 \pm 4.4) \times 10^{-5}\) & \((1.2 \pm 1.1) \times 10^{-5}\) \\
		\(\Loss_{\mse}\) training data & \((7.9 \pm 4.6) \times 10^{-6}\) & \((5.3 \pm 0.5) \times 10^{-6}\) \\
		\(\Loss_{\mse}\) validation data & \((1.7 \pm 4.4) \times 10^{-5}\) & \((3.8 \pm 2.4) \times 10^{-6}\) \\
		number of epochs & \(74 \pm 19\) & \(39 \pm 14\) \\
		training time (in \(s\)) & \(1307 \pm 346\) & \(2075\pm 829\) \\
		\textbf{speedup factor} & \textbf{1.6} & \\
		\midrule
	\end{tabular}
	\caption{
		Summary of the learning statistics for the neural network compression approaches \(\Upsilon^{\pc}_{\NN}\) and \(\Phi^{\pc}_{\NN}\) for the four layer network \(k_d \rightarrow 10 \rightarrow 10 \rightarrow 20 \rightarrow 1\). 
		Results are averaged over ten runs and additionally include the standard deviation. 
	}
	\label{tab:NN_Comparison}
\end{table}

To complement this qualitative assessment, \cref{tab:NN_Comparison} reports quantitative performance measures of the training process for both approaches. 
The training and validation losses of both approaches are around the same order of magnitude. 
However, a direct comparison of the loss values is only partially meaningful due to the differing symmetry notions and the fact that, for instance, the mean squared error is evaluated on different computational domains. 
For this reason, additional illustrations provided below aim to shed further light on the approximation capabilities of the respective models. 

In terms of training duration, the \(\Upsilon^{\pc}_{\NN}\) approach typically requires a larger number of epochs, it nevertheless achieves a significantly lower overall training time. 
In the present configuration, the average training time is approximately \num{22} minutes for \(\Upsilon^{\pc}_{\NN}\) and \num{35} minutes for \(\Phi^{\pc}_{\NN}\), corresponding to a speedup factor of about \num{1.6}. 
We emphasise that this speedup should be regarded as conservative. 
Across a range of additional architectures and penalty parameter choices, we observed speedup factors between approximately \(1.4\) and \(6\). 
The precise value depends on the chosen hyperparameters and training configuration. 
\smallskip

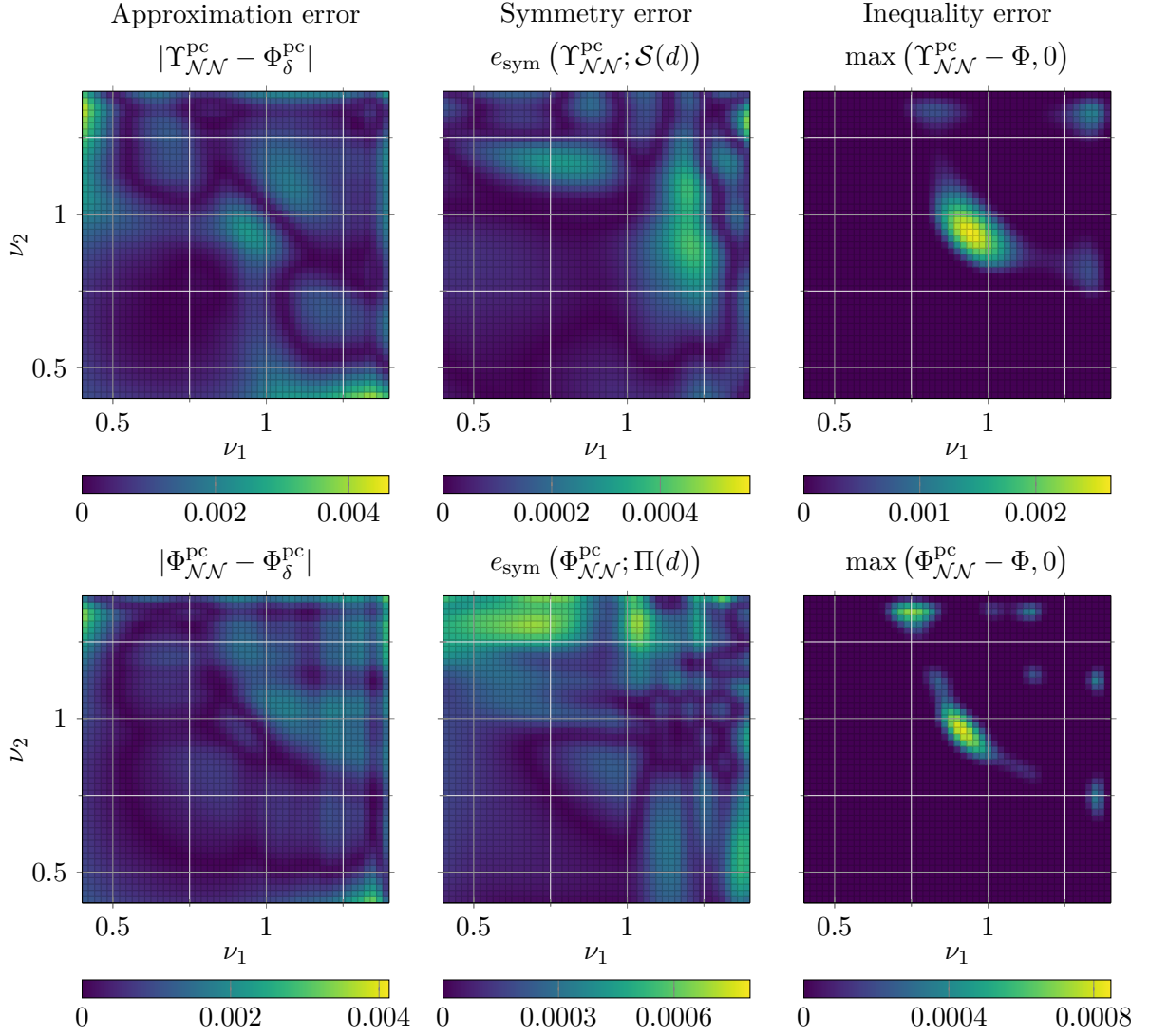
\begin{figure}
	\centering	
	\begin{tikzpicture}
	
	\pgfmathsetmacro{\widthscale}{0.31}		
	\pgfmathsetmacro{\xshift}{0.75} 
	\pgfmathsetmacro{\yshift}{-2.75} 
	
	\begin{axis}[
		width=\widthscale\textwidth,
		scale only axis,
		grid=both, 
		tick align=center,
		minor grid style={black!10},
		major grid style={black!40},
		axis background/.style={fill=white},
		view={0}{90},
		xmajorgrids, ymajorgrids,
		xlabel={\(\nu_1\)},	
		ylabel={\(\nu_2\)},
		xlabel style={yshift=1.mm},
		ylabel style={yshift=-2.mm},
		xtick={{0,0.5,1,1.5}}, 
		ytick={{0,0.5,1,1.5}}, 
		minor x tick num=1,
		minor y tick num=1,
		colormap/viridis,
		mesh/ordering=x varies,
		colorbar horizontal,
		axis equal image,
		name=plotA,
		colorbar style={
			height=2.5mm,
			xticklabel style={
				/pgf/number format/fixed,
				/pgf/number format/precision=4
			},
			yshift=2mm,
			scaled ticks=false,	
		}
		]
		\addplot3[surf, shader=faceted] file {figures/data/PhipcNNsliceSTVK_err_3d_Ball.dat};
	\end{axis}
	\node[anchor=south, yshift=7.5mm] at (plotA.north) {Approximation error};
	\node[anchor=south, yshift=1mm] at (plotA.north) {\(\lvert \Upsilon^{\pc}_{\NN} - \Phi^{\pc}_{\delta}\rvert\)};
	
	\begin{axis}[
		width=\widthscale\textwidth,
		scale only axis,
		grid=both,
		tick align=center, 
		minor grid style={black!10},
		major grid style={black!40},
		axis background/.style={fill=white},
		view={0}{90},
		xmajorgrids, ymajorgrids,
		xlabel={\(\nu_1\)},	
		xlabel style={yshift=1.mm},
		ylabel style={yshift=-2.mm},
		xtick={{0,0.5,1,1.5}}, 
		ytick={{0,0.5,1,1.5}}, 
		minor x tick num=1,
		minor y tick num=1,
		yticklabels={},
		colormap/viridis,
		mesh/ordering=x varies,
		colorbar horizontal,
		axis equal image,
		at={(plotA.south east)},
		xshift=\xshift cm,
		name=plotB,
		colorbar style={
			height=2.5mm,
			xticklabel style={
				/pgf/number format/fixed,
				/pgf/number format/precision=4
			},
			yshift=2mm,
			scaled ticks=false,
		}
		]
		\addplot3[surf, shader=faceted] file {figures/data/PhipcNNsliceSTVK_errsym_3d_Ball.dat};
	\end{axis}
	\node[anchor=south, yshift=7.5mm] at (plotB.north) {Symmetry error};
	\node[anchor=south, yshift=1mm] at (plotB.north) {\(e_{\sym}\left(\Upsilon^{\pc}_{\NN}; \Sd\right)\)};

	\begin{axis}[
		width=\widthscale\textwidth,
		scale only axis,
		grid=both, 
		tick align=center,
		minor grid style={black!10},
		major grid style={black!40},
		axis background/.style={fill=white},
		view={0}{90},
		xmajorgrids, ymajorgrids,
		xlabel={\(\nu_1\)},	
		xlabel style={yshift=1.mm},
		ylabel style={yshift=-2.mm},
		xtick={{0,0.5,1,1.5}}, 
		ytick={{0,0.5,1,1.5}}, 
		minor x tick num=1,
		minor y tick num=1, 
		yticklabels={},
		colormap/viridis,
		mesh/ordering=x varies,
		colorbar horizontal,
		axis equal image,
		at={(plotB.south east)},
		xshift=\xshift cm,
		name=plotC,
		colorbar style={
			height=2.5mm,
			xticklabel style={
				/pgf/number format/fixed,
				/pgf/number format/precision=4
			},
			yshift=2mm,
			scaled ticks=false,
		}
		]
		\addplot3[surf, shader=faceted] file {figures/data/PhipcNNsliceSTVK_errineq_3d_Ball.dat};
	\end{axis}
	\node[anchor=south, yshift=7.5mm] at (plotC.north) {Inequality error};
	\node[anchor=south, yshift=1mm] at (plotC.north) {\(\max\left(\Upsilon^{\pc}_{\NN} - \Phi, 0\right)\)};
	
	\begin{axis}[
		width=\widthscale\textwidth,
		scale only axis,
		grid=both, 
		tick align=center,
		minor grid style={black!10},
		major grid style={black!40},
		axis background/.style={fill=white},
		view={0}{90},
		xmajorgrids, ymajorgrids,
		xlabel={\(\nu_1\)},	
		ylabel={\(\nu_2\)},
		xlabel style={yshift=1.mm},
		ylabel style={yshift=-2.mm},
		xtick={{0,0.5,1,1.5}}, 
		ytick={{0,0.5,1,1.5}}, 
		minor x tick num=1,
		minor y tick num=1,
		colormap/viridis,
		mesh/ordering=x varies,
		colorbar horizontal,
		axis equal image,
		at={(plotA.south west)},
		anchor=north west,
		yshift=\yshift cm,
		name=plotD,
		colorbar style={
			height=2.5mm,
			xtick={0,0.002,0.004},
			xticklabel style={
				/pgf/number format/fixed,
				/pgf/number format/precision=4
			},
			yshift=2mm,
			scaled ticks=false,	
		}
		]
		\addplot3[surf, shader=faceted] file {figures/data/PhipcNNsliceSTVK_err_3d_Iso.dat};
	\end{axis}
	\node[anchor=south, yshift=1mm] at (plotD.north) {\(\lvert \Phi^{\pc}_{\NN} - \Phi^{\pc}_{\delta}\rvert\)};
	
	\begin{axis}[
		width=\widthscale\textwidth,
		scale only axis,
		grid=both,
		tick align=center, 
		minor grid style={black!10},
		major grid style={black!40},
		axis background/.style={fill=white},
		view={0}{90},
		xmajorgrids, ymajorgrids,
		xlabel={\(\nu_1\)},	
		xlabel style={yshift=1.mm},
		ylabel style={yshift=-2.mm},
		xtick={{0,0.5,1,1.5}}, 
		ytick={{0,0.5,1,1.5}}, 
		minor x tick num=1,
		minor y tick num=1,
		yticklabels={},
		colormap/viridis,
		mesh/ordering=x varies,
		colorbar horizontal,
		axis equal image,
		at={(plotD.south east)},
		xshift=\xshift cm,
		name=plotE,
		colorbar style={
			height=2.5mm,
			xtick={0,0.0003,0.0006},
			xticklabel style={
				/pgf/number format/fixed,
				/pgf/number format/precision=4
			},
			yshift=2mm,
			scaled ticks=false,
		}
		]
		\addplot3[surf, shader=faceted] file {figures/data/PhipcNNsliceSTVK_errsym_3d_Iso.dat};
	\end{axis}
	\node[anchor=south, yshift=1mm] at (plotE.north) {\(e_{\sym}\left(\Phi^{\pc}_{\NN}; \Pid\right)\)};

	\begin{axis}[
		width=\widthscale\textwidth,
		scale only axis,
		grid=both, 
		tick align=center,
		minor grid style={black!10},
		major grid style={black!40},
		axis background/.style={fill=white},
		view={0}{90},
		xmajorgrids, ymajorgrids,
		xlabel={\(\nu_1\)},	
		xlabel style={yshift=1.mm},
		ylabel style={yshift=-2.mm},
		xtick={{0,0.5,1,1.5}}, 
		ytick={{0,0.5,1,1.5}}, 
		minor x tick num=1,
		minor y tick num=1, 
		yticklabels={},
		colormap/viridis,
		mesh/ordering=x varies,
		colorbar horizontal,
		axis equal image,
		at={(plotE.south east)},
		xshift=\xshift cm,
		name=plotF,
		colorbar style={
			height=2.5mm,
			xtick={0,0.0004,0.0008},
			xticklabel style={
				/pgf/number format/fixed,
				/pgf/number format/precision=4
			},
			yshift=2mm,
			scaled ticks=false,
		}
		]
		\addplot3[surf, shader=faceted] file {figures/data/PhipcNNsliceSTVK_errineq_3d_Iso.dat};
	\end{axis}
	\node[anchor=south, yshift=1mm] at (plotF.north) {\(\max\left(\Phi^{\pc}_{\NN} - \Phi, 0\right)\)};

\end{tikzpicture}
	\caption{
		Pointwise approximation error (left column), pointwise symmetry error (middle column) and pointwise inequality error (right column) of the neural network representations \(\Upsilon^{\pc}_{\NN}\) (top row) and \(\Phi^{\pc}_{\NN}\) (bottom row) on the \((\nu_1, \nu_2, 1)\) cross section for \(\nu_1, \nu_2 \in [0.4, 1.4]^{2}\).
		Network outputs are averaged over ten network realisations. 
	}
	\label{fig:NNpropertyerrorsSTVK3D}
\end{figure}

\cref{fig:NNpropertyerrorsSTVK3D} shows the individual approximation and property errors of the neural network compressions.
The top row corresponds to the compression via \(\Upsilon^{\pc}_{\NN}\), while the bottom row shows the results for \(\Phi^{\pc}_{\NN}\).
The left column depicts the pointwise approximation error, the middle column the symmetry error, and the right column the inequality error, respectively.

As a local measure of non-symmetry of the neural network compression, we consider the pointwise deviation from the symmetrised network output.
For a network \(\hat{\Gamma}^{\pc}_{\NN} \in \{\Upsilon^{\pc}_{\NN}, \Phi^{\pc}_{\NN}\}\) with associated symmetry group \(\gensym \in \{\Sd, \Pid\}\), the symmetry error \(e_{\sym}(\hat{\Gamma}^{\pc}_{\NN}; \gensym) \colon \R^{3} \to \R\) describes the pointwise distance to the symmetrised network output and is defined by
\begin{equation*}
	e_{\sym}(\hat{\Gamma}^{\pc}_{\NN}; \gensym)(\hat{\gamma}) = \hat{\Gamma}^{\pc}_{\NN}(\hat{\gamma}) - \frac{1}{\lvert \gensym \rvert} \sum_{S \in \gensym}\hat{\Gamma}^{\pc}_{\NN}(S \, \hat{\gamma}).
\end{equation*} 

The pointwise approximation errors indicate comparable accuracy for both compression approaches. 
Moreover, both methods exhibit symmetry and inequality errors of similar magnitude, significantly lower than the pure pointwise approximation error, confirming effective incorporation of the structural properties through the penalty formulation.

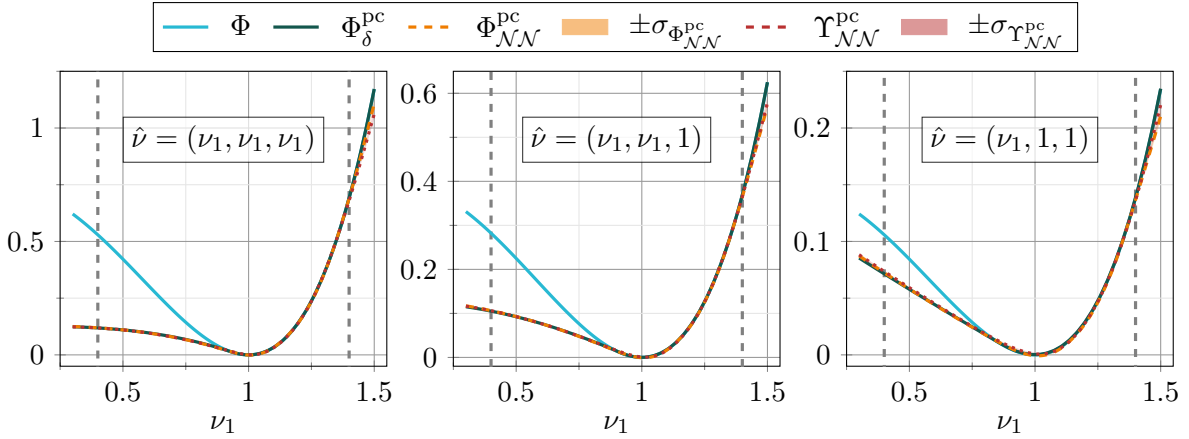
\begin{figure}
	\centering	
		\begin{tikzpicture}
	
	\pgfmathsetmacro{\widthscale}{0.27}	
	\pgfmathsetmacro{\xshift}{0.75cm} 
	\pgfmathsetmacro{\height}{3.9cm}
	
	\colorlet{colPhi}{unia-lightblue}
	\colorlet{colPhipcdelta}{unia-green}
	\colorlet{colPhipcNNIso}{unia-orange}
	\colorlet{colPhipcNNBall}{unia-red}
	
	\begin{axis}[
		width=\widthscale\textwidth,
		height=\height,
		scale only axis,		
		xlabel={\(\nu_1\)},
		grid=both, 
		tick align=center,
		tick pos=left,
		minor grid style={black!10},
		major grid style={black!40},
		axis background/.style={fill=white},
		xmajorgrids, ymajorgrids,
		minor x tick num=1,
		minor y tick num=1,
		xmin=0.25, xmax=1.55,
		ymin=-0.05, ymax=1.25,
		legend cell align={left},
		legend style={
			fill opacity=1,
			draw opacity=1,
			text opacity=1,
			at={(1.7,1.15)},
			anchor=center,
			legend columns=6,
			column sep=2mm,
		},
		at={(0, 0)},
		name=plotA1dcross,	
		]

		\addlegendimage{colPhi, very thick}
		\addlegendentry{\({\Phi}\)}
		
		\addlegendimage{colPhipcdelta, very thick}
		\addlegendentry{\({\Phi^{\pc}_{\delta}}\)}
		
		\addlegendimage{colPhipcNNIso, dashed, very thick}
		\addlegendentry{\({\Phi^{\pc}_{\NN}}\)}
		
		\addlegendimage{area legend, draw=colPhipcNNIso, fill=colPhipcNNIso, opacity=0.5}
		\addlegendentry{\(\pm \sigma_{\Phi^{\pc}_{\NN}}\)}
		
		\addlegendimage{colPhipcNNBall, dashed, very thick}
		\addlegendentry{\({\Upsilon^{\pc}_{\NN}}\)}
		
		\addlegendimage{area legend, draw=colPhipcNNBall, fill=colPhipcNNBall, opacity=0.5}
		\addlegendentry{\(\pm \sigma_{\Upsilon^{\pc}_{\NN}}\)}
		
		\addplot[very thick, dashed, gray]
		coordinates {
			(0.4,-1)
			(0.4, 2)
		};
		
		\addplot[very thick, dashed, gray]
		coordinates {
			(1.4,-1)
			(1.4, 2)
		};
	
		\addplot[very thick, colPhi]
		table[
		x index=0,
		y index=1
		] {figures/data/Phi_tri.dat};

		\addplot[very thick, colPhipcdelta]
		table[
		x index=0,
		y index=1
		] {figures/data/Phipc_delta_tri.dat};
		
		\addplot[name path=lowerNNIso, draw=none]
		table[x=x,y=lower] {figures/data/Phipc_NN_Iso_tri.dat};
		
		\addplot[name path=upperNNIso, draw=none]
		table[x=x,y=upper] {figures/data/Phipc_NN_Iso_tri.dat};
		
		\addplot[fill opacity=0.2, fill=colPhipcNNIso]
		fill between[of=lowerNNIso and upperNNIso];
		
		\addplot[very thick, colPhipcNNIso, dashed]
		table[x=x,y=mean] {figures/data/Phipc_NN_Iso_tri.dat};
		
		\addplot[name path=lowerNNBall, draw=none]
		table[x=x,y=lower] {figures/data/Phipc_NN_Ball_tri.dat};
		
		\addplot[name path=upperNNBall, draw=none]
		table[x=x,y=upper] {figures/data/Phipc_NN_Ball_tri.dat};
		
		\addplot[fill opacity=0.2, fill=colPhipcNNBall]
		fill between[of=lowerNNBall and upperNNBall];
		
		\addplot[very thick, colPhipcNNBall, dotted]
		table[x=x,y=mean] {figures/data/Phipc_NN_Ball_tri.dat};
				
	\end{axis}
	
	\begin{axis}[
		width=\widthscale\textwidth,
		height=\height,
		scale only axis,
		xlabel={\(\nu_1\)},
				grid=both, 
		tick align=center,
		tick pos=left,
		minor grid style={black!10},
		major grid style={black!40},
		axis background/.style={fill=white},
		xmajorgrids, ymajorgrids,
		minor x tick num=1,
		minor y tick num=1,
		xmin=0.25, xmax=1.55,
		ymin=-0.02, ymax=0.65,
		at={(0.325\textwidth, 0)},
		name=plotB1dcross,
		]
		
		\addplot[very thick, dashed, gray]
		coordinates {
			(0.4,-1)
			(0.4, 2)
		};
		
		\addplot[very thick, dashed, gray]
		coordinates {
			(1.4,-1)
			(1.4, 2)
		};
	
		\addplot[very thick, colPhi]
		table[
		x index=0,
		y index=1
		] {figures/data/Phi_bi.dat};

		\addplot[very thick, colPhipcdelta]
		table[
		x index=0,
		y index=1
		] {figures/data/Phipc_delta_bi.dat};

		\addplot[name path=lowerNNIso, draw=none]
		table[x=x,y=lower] {figures/data/Phipc_NN_Iso_bi.dat};
		
		\addplot[name path=upperNNIso, draw=none]
		table[x=x,y=upper] {figures/data/Phipc_NN_Iso_bi.dat};
		
		\addplot[fill opacity=0.2, fill=colPhipcNNIso]
		fill between[of=lowerNNIso and upperNNIso];
		
		\addplot[very thick, colPhipcNNIso, dashed]
		table[x=x,y=mean] {figures/data/Phipc_NN_Iso_bi.dat};
		
		\addplot[name path=lowerNNBall, draw=none]
		table[x=x,y=lower] {figures/data/Phipc_NN_Ball_bi.dat};
		
		\addplot[name path=upperNNBall, draw=none]
		table[x=x,y=upper] {figures/data/Phipc_NN_Ball_bi.dat};
		
		\addplot[fill opacity=0.2, fill=colPhipcNNBall]
		fill between[of=lowerNNBall and upperNNBall];
		
		\addplot[very thick, colPhipcNNBall, dotted]
		table[x=x,y=mean] {figures/data/Phipc_NN_Ball_bi.dat};
		
	\end{axis}
	
	\begin{axis}[
		width=\widthscale\textwidth,
		height=\height,
		scale only axis,
		xlabel={\(\nu_1\)},
		grid=both, 
		tick align=center,
		tick pos=left,
		minor grid style={black!10},
		major grid style={black!40},
		axis background/.style={fill=white},
		xmajorgrids, ymajorgrids,
		minor x tick num=1,
		minor y tick num=1,
		xmin=0.25, xmax=1.55,
		ymin=-0.01, ymax=0.25,
		at={(0.65\textwidth, 0)},
		name=plotC1dcross,
		]
		
		\addplot[very thick, dashed, gray]
		coordinates {
			(0.4,-1)
			(0.4, 2)
		};
		
		\addplot[very thick, dashed, gray]
		coordinates {
			(1.4,-1)
			(1.4, 2)
		};
		
		\addplot[very thick, colPhi]
		table[
		x index=0,
		y index=1
		] {figures/data/Phi_uni.dat};

		\addplot[very thick, colPhipcdelta]
		table[
		x index=0,
		y index=1
		] {figures/data/Phipc_delta_uni.dat};

		\addplot[name path=lowerNNIso, draw=none]
		table[x=x,y=lower] {figures/data/Phipc_NN_Iso_uni.dat};
		
		\addplot[name path=upperNNIso, draw=none]
		table[x=x,y=upper] {figures/data/Phipc_NN_Iso_uni.dat};
		
		\addplot[fill opacity=0.2, fill=colPhipcNNIso]
		fill between[of=lowerNNIso and upperNNIso];
		
		\addplot[very thick, colPhipcNNIso, dashed]
		table[x=x,y=mean] {figures/data/Phipc_NN_Iso_uni.dat};
		
		\addplot[name path=lowerNNBall, draw=none]
		table[x=x,y=lower] {figures/data/Phipc_NN_Ball_uni.dat};
		
		\addplot[name path=upperNNBall, draw=none]
		table[x=x,y=upper] {figures/data/Phipc_NN_Ball_uni.dat};
		
		\addplot[fill opacity=0.2, fill=colPhipcNNBall]
		fill between[of=lowerNNBall and upperNNBall];
		
		\addplot[very thick, colPhipcNNBall, dotted]
		table[x=x,y=mean] {figures/data/Phipc_NN_Ball_uni.dat};
		
	\end{axis}
	
	\node[draw, 
	fill=white!30, 
	draw opacity=0.8, 
	fill opacity=0.7, 
	text opacity=1, 
	inner sep=3pt, 
	anchor=north, 
	yshift=-0.6cm, 
	] at (plotA1dcross.north) {\(\hat{\nu} = (\nu_1, \nu_1, \nu_1)\)};
	
	\node[draw, 
	fill=white!30, 
	draw opacity=0.8, 
	fill opacity=0.7, 
	text opacity=1, 
	inner sep=3pt, 
	anchor=north, 
	yshift=-0.6cm, 
	] at (plotB1dcross.north) {\(\hat{\nu} = (\nu_1, \nu_1, 1)\)};
	
	\node[draw, 
	fill=white!30, 
	draw opacity=0.8, 
	fill opacity=0.7, 
	text opacity=1, 
	inner sep=3pt, 
	anchor=north, 
	yshift=-0.6cm, 
	] at (plotC1dcross.north) {\(\hat{\nu} = (\nu_1, 1, 1)\)};

\end{tikzpicture}
	\caption{
		One-dimensional cross sections for the three-dimensional Saint~Venant--Kirchhoff example. 
		The two compression strategies \(\Upsilon^{\pc}_{\NN}\) and \(\Phi^{\pc}_{\NN}\) are the results of ten averaged network realisations with indicated (almost vanishing) standard deviations. 
		Illustrated is the domain \(\nu_1 \in [0.3, 1.5]\) where the evaluation is performed on 200 equidistant points. 
		Additionally, the boundary of the learning domain (\(\nu_1 \in [0.4, 1.4]\)) is indicated.
	}
	\label{fig:1DcrosssectionsSTVK3D}
\end{figure}

\cref{fig:1DcrosssectionsSTVK3D} illustrates three one-dimensional cross sections of the three-dimensional Saint~Venant--Kirchhoff example for the neural network compressions \(\Upsilon^{\pc}_{\NN}\) and \(\Phi^{\pc}_{\NN}\).
Both approaches exhibit a reliable approximation of the reference envelope, with no discernible difference in accuracy along the selected slices.
The variability across independent training runs is negligible, as indicated by the virtually vanishing standard deviations \(\sigma_{\Phi^{\pc}_{\NN}}\) and \(\sigma_{\Upsilon^{\pc}_{\NN}}\), justifying the presentation of averaged network outputs in the preceding figures.
Both compression approaches demonstrate stable and accurate interpolation within the training domain, while maintaining acceptable accuracy in regions outside the training domain.

Overall, both approaches yield accurate and structure-preserving approximations of the polyconvex envelope, with no significant differences in accuracy and in the enforcement of the required physical properties. 
At the same time, the neural network reduces the \num{1687500} (resp.~\num{421875}) data points for \(\Phi^{\pc}_{\NN}\) (resp.~\(\Upsilon^{\pc}_{\NN}\)) to a representation with only \num{648} trainable parameters, which corresponds, in terms of storage capacity, to representing the envelope on a \(9 \times 9 \times 9\) lattice and highlights the huge compression potential for precomputed polyconvex envelopes.  
The Ball-criterion-based formulation is characterised by a consistently reduced training time and, consequently, a significant decrease in computational cost.

\section{Conclusion}
The proposed structure-preserving neural network framework provides an efficient and reliable approach for the compression of polyconvex envelopes of isotropic functions based on the classical sufficient criterion for polyconvexity of \cite{Bal76}. 
Although this compression approach generally provide only a lower bound on the polyconvex envelope, the proposed architecture proved successful for the physically relevant determinant-constrained energy density considered in this work, whose polyconvex envelope is not available in analytical form. 
Substantial reductions in computational cost can be achieved through the reduction from the signed singular value space to the singular value space while preserving the required structural properties and maintaining the accuracy of state-of-the-art approaches.
These computational savings appear particularly attractive in parameter-dependent settings, where they may facilitate concurrent relaxation in boundary value problem simulations, such as damage models, a direction that could be pursued in future work.

\section*{Acknowledgement}
The authors would like to express their gratitude for the fruitful scientific discussions with Daniel Peterseim and Malte A.~Peter. 

\printbibliography

\end{document}